\documentclass{amsart}
\usepackage[utf8]{inputenc}
\usepackage{latexsym,amssymb,amsmath,amsthm,amscd,graphicx}
\usepackage{amsfonts}
\usepackage{xcolor}
\usepackage[english]{babel}
\usepackage{color}
\definecolor{dblue}{rgb}{0,0,0.45}
\usepackage[a4paper,bindingoffset=0.2in,%
            left=1in,right=1in,top=1in,bottom=1in,%
            footskip=.25in]{geometry}
\newtheorem{theorem}{Theorem}[section]
\newtheorem{lemma}[theorem]{Lemma}

\newtheorem{cor}[theorem]{Corollary}
\begin{document}

\title[On geometric properties of Morrey spaces]{On geometric properties of Morrey spaces}

\author[H.~Gunawan]{Hendra Gunawan}
\address{Faculty of Mathematics and Natural Sciences, Bandung Institute of Technology, Bandung 40132, Indonesia}
\email{hgunawan@math.itb.ac.id}

\author[D.I.~Hakim]{Denny I. Hakim}
\address{Faculty of Mathematics and Natural Sciences, Bandung Institute of Technology, Bandung 40132, Indonesia}
\email{dhakim@math.itb.ac.id}

\author[A.S.~Putri]{Arini S. Putri}
\address{Faculty of Mathematics and Natural Sciences, Bandung Institute of Technology, Bandung 40132, Indonesia}
\email{arini@students.itb.ac.id}



\begin{abstract}
    In this article, we show constructively that Morrey spaces are not uniformly non-$\ell^1_n$
    for any $n\ge 2$. This result is sharper than those previously obtained in \cite{GKSS, MG},
    which show that Morrey spaces are not uniformly non-square and also not uniformly non-octahedral.
    We also discuss the $n$-th James constant $C_{{\rm J}}^{(n)}(X)$ and the $n$-th Von
    Neumann-Jordan constant $C_{{\rm NJ}}^{(n)}(X)$ for a Banach space $X$, and obtain
    that both constants for any Morrey space $\mathcal{M}^p_q(\mathbb{R}^d)$ with
    $1\le p<q<\infty$ are equal to $n$.
\end{abstract}

\maketitle

\section{Introduction}

For $1\leq p\leq q<\infty$, the {\em Morrey space} $\mathcal{M}^p_q=\mathcal{M}^p_q(\mathbb{R}^d)$
is the set of all measurable functions $f$ such that
\begin{equation*}
    \|f\|_{\mathcal{M}^p_q} := \sup_{{a\in \mathbb{R}^d},{R>0}}|B(a,R)|^{\frac{1}{q}-\frac{1}{p}}
    \biggl(\int\limits_{B(a,R)} |f(y)|^p dy\biggr)^{\frac{1}{p}}<\infty,
\end{equation*}
where $|B(a,R)|$ denotes the Lebesgue measure of the open ball $B(a,R)$ in $\mathbb{R}^d$, with center $a$ and
radius $R$. Morrey spaces are Banach spaces (see, e.g., \cite{sawano}). For $p=q$, the space $\mathcal{M}^q_q$ is
identical with the space $L^q=L^q(\mathbb{R}^d)$, the space of $q$-th power integrable functions on $\mathbb{R}^d$.

In \cite{GKSS}, three geometric constants have been computed for Morrey spaces. The first two constants, namely
{\em Von Neumann--Jordan constant} and {\em James constant}, are closely related to the notion of uniformly
non-squareness of (the unit ball of) a Banach space \cite{james, jimenez, kato}. For a Banach space $(X,\|\cdot\|_X)$
in general, the constants are defined by
\begin{equation*}
    C_{\rm NJ}(X) := \sup\biggl\{\frac{\|x+y\|_X^2+\|x-y\|_X^2}{2(\|x\|_X^2+\|y\|_X^2)} : x,y \in X\setminus \{0\} \biggr\},
\end{equation*}
and
\begin{equation*}
    C_{\rm J}(X) := \sup\bigl\{\min\{\|x+y\|_X,\|x-y\|_X\} : x,y\in S_X \bigr\},
\end{equation*}
respectively. Here $S_X:=\{x\in X : \|x\|_X=1\}$ denotes the unit sphere in $X$. A few basic facts about these constants are:
\begin{itemize}
    \item $1\leq C_{\rm NJ}(X)\leq 2$ and $C_{\rm NJ}(X)=1$ if and only if $X$ is a Hilbert space \cite{jordan1935inner}.
    \item $\sqrt{2}\leq C_{\rm J}(X)\leq 2$ and $C_{\rm J}(X)=\sqrt{2}$ if $X$ is a Hilbert space \cite{gao1990geometry}.
\end{itemize}
Note also that, for $1\le p\le \infty$, we have $C_{\rm NJ}(L^p)=\max\{2^{\frac{2}{p}-1},2^{1-\frac{2}{p}}\}$
and $C_{\rm J}(L^p)=\max\{2^{\frac{1}{p}},2^{1-\frac{1}{p}}\}$ \cite{clarkson, gao1990geometry}.

As for Morrey spaces, we have the following results.

\medskip

\begin{theorem}{\rm \cite{GKSS}}\label{Morrey}
For $1\le p<q<\infty$, we have $C_{\rm NJ}(\mathcal{M}^p_q)=C_{\rm J}(\mathcal{M}^p_q)=2$.
\end{theorem}

\medskip

This theorem tells us that Morrey spaces are not {\em uniformly non-square} (since a Banach space is uniformly non-square
if and only if $C_{\rm J}(X)<2$ or, equivalently, $C_{\rm NJ}(X)<2$). In \cite{MG}, we are also told that Morrey spaces
are not {\em uniformly non-octahedral}, that is, there does not exist a $\delta>0$ such that
\[
\min \|f\pm g\pm h\|_{\mathcal{M}^p_q} \le 3(1-\delta)
\]
for every $f,g,h\in \mathcal{M}^p_q$ with $\|f\|_{\mathcal{M}^p_q}=\|g\|_{\mathcal{M}^p_q}=\|h\|_{\mathcal{M}^p_q}=1$.
Here the minimum is taken over all choices of signs in the expression $f\pm g\pm h$. (A Banach space $(X,\|\cdot\|_X)$
is uniformly non-octahedral if there exists a $\delta>0$ such that
\[
\min \|x\pm y\pm z\|_X \le 3(1-\delta)
\]
for every $x,y,z\in S_X$.) Precisely, we have the following theorem.

\medskip

\begin{theorem}{\rm \cite{MG}}\label{Octahedral}
Let $1\leq p< q<\infty$. Then, for every $\delta>0$, there exist $f,g,h\in \mathcal{M}^p_q$ (depending
on $\delta$) with $\|f\|_{\mathcal{M}^p_q}=\|g\|_{\mathcal{M}^p_q}=\|h\|_{\mathcal{M}^p_q}=1$ such that
\[
\|f\pm g\pm h\|_{\mathcal{M}^p_q}> 3(1-\delta)
\]
for all choices of signs.
\end{theorem}

\medskip

In this article, we shall show constructively that Morrey spaces are not {\em uniformly non}-$\ell^1_n$. The result is not only more
general than the previous ones, but also sharper than knowing that Morrey spaces are neither uniformly non-square nor
uniformly non-octahedral (for if $X$ is not uniformly non-$\ell^1_n$ for $n\ge3$, then $X$ is not uniformly non-$\ell^1_{n-1}$).
In addition, given a Banach space $X$, we shall discuss the $n$-th Von Neumann-Jordan constant $C_{{\rm NJ}}^{(n)}(X)$ and
the $n$-th James constant $C_{{\rm J}}^{(n)}(X)$ for $n\ge2$. These two constants were studied in \cite{KTH} and \cite{MNPZ},
respectively. We show that for any Morrey space $\mathcal{M}^p_q$ with $1\le p<q<\infty$ both constants are equal to $n$.
We also indicate that Morrey spaces are not {\em uniformly $n$-convex} for $n\ge2$.

\section{$\mathcal{M}^p_q$ is not uniformly non-$\ell^1_n$}

Before we present our main theorems, we prove several lemmas. We assume that the readers know well how to compute
the integral of a radial function over a ball centered at 0 using polar coordinates. Unless otherwise stated,
we assume that $1\le p<q<\infty$.

\medskip

\begin{lemma}\label{le:140320-1}
Let $f(x):=|x|^{-d/q}$. Then $f\in \mathcal{M}^p_q$ with
\[
\|f\|_{\mathcal{M}^p_q}=\Bigl(\frac{\omega_{d-1}}{d}\Bigr)^{\frac{1}{q}} \Bigl(\frac{q}{q-p}\Bigr)^{\frac{1}{p}},
\]
where $\omega_{d-1}$ denotes the `area' of the unit sphere in $\mathbb{R}^d$.
\end{lemma}

\medskip

\begin{proof}
For each $r>0$, one may compute that
\[
|B(0,r)|^{1/q-1/p}\Bigl(\int_{B(0,r)} |x|^{-dp/q}dx\Bigr)^{1/p}=\Bigl(\frac{\omega_{d-1}}{d}\Bigr)^{1/q}
\Bigl(\frac{q}{q-p}\Bigr)^{1/p},
\]
which is independent of $r$. Since the integral of $f$ over $B(a,r)$ will be less than that over $B(0,r)$
for every $a\in\mathbb{R}^d$, we conclude that $\|f\|_{\mathcal{M}^p_q}=
\bigl(\frac{\omega_{d-1}}{d}\bigr)^{\frac{1}{q}} \bigl(\frac{q}{q-p}\bigr)^{\frac{1}{p}}$, as claimed.
\end{proof}

\medskip

\begin{lemma}\label{le:140320-2}
Let $f(x):=|x|^{-d/q}$ and $R>1$. Then, for any $c_1,c_2>0$, we have
\begin{align}\label{eq:140320-3}
|B_{c_1R}|^{1/q-1/p}\Bigl(\int_{\{x:c_1<|x|<c_1R\}} |f(x)|^p dx\Bigr)^{1/p}=|B_{c_2R}|^{1/q-1/p}
\Bigl(\int_{\{x:c_2<|x|<c_2R\}}|f(x)|^p dx\Bigr)^{1/p},
\end{align}
where $B_{c_1R}$ and $B_{c_2R}$ are balls centered at $0$ with radii $c_1R$ and $c_2R$.
\end{lemma}

\medskip

\begin{proof}
It suffices to prove that \eqref{eq:140320-3} holds for arbitrary $c_1>0$ and $c_2=1$. But this is
immediate by the change of variable $x=c_1x'$.
\end{proof}

As a consequence of the above lemma, we have the following corollary, which is an important key to
our main theorems.

\medskip

\begin{cor}\label{cor:140320-1}
Let $f(x):=|x|^{-d/q}$. For $\varepsilon\in(0,1)$ and $k\in\mathbb{N}$, put $f_{\varepsilon,k}:=
f\chi_{\{x:\varepsilon^{k+1}<|x|<\varepsilon^{k}\}}$. Then $f_{\varepsilon,k}\in \mathcal{M}^p_q$
with
\begin{align}\label{eq:140320-2}
\|f_{\varepsilon,k}\|_{\mathcal{M}^p_q}\ge (1-\varepsilon^{d-dp/q})^{1/p}\|f\|_{\mathcal{M}^p_q}.
\end{align}
\end{cor}

\medskip

\begin{proof}
In view of Lemma \ref{le:140320-2}, it suffices to prove that
\[
\|f_{\varepsilon,0}\|_{\mathcal{M}^p_q}\ge (1-\varepsilon^{d-dp/q})^{1/p}\|f\|_{\mathcal{M}^p_q}.
\]
We observe that
\begin{align*}
\|f_{\varepsilon,0}\|_{\mathcal{M}^p_q}
\ge |B(0,1)|^{\frac{1}{q}-\frac1p}\Bigl(\int_{\{x:\varepsilon<|x|<1\}} |f(x)|^pdx\Bigr)^{1/p}=
\Bigl(\frac{\omega_{d-1}}{d}\Bigr)^{\frac{1}{q}} \Bigl(\frac{q}{q-p}\Bigr)^{\frac{1}{p}}(1-\varepsilon^{d-dp/q})^{1/p}.
\end{align*}
Hence, by Lemma \ref{le:140320-1}, the desired inequality follows.
\end{proof}

We are now ready to state our main results. Our first theorem is the following.

\medskip

\begin{theorem}\label{unl1n}
For $1\leq p< q<\infty$, the Morrey space $\mathcal{M}^p_q$ is not uniformly non-$\ell^1_n$ for any $n\ge 2$,
that is, for every $\delta\in (0,1)$, there exist $f_1,f_2\dots, f_n\in \mathcal{M}^p_q$
(depending on $\delta$) with $\|f_i\|_{\mathcal{M}^p_q}=1$ for $i=1,2,\dots,n$, such that
\[
\|f_1\pm f_2\pm \cdots \pm f_n\|_{\mathcal{M}^p_q}> n(1-\delta)
\]
for all choices of signs.
\end{theorem}

\medskip

\begin{proof}
To understand the idea of the proof, let us first ilustrate how the proof goes for $n=3$. Given $\delta \in (0,1)$,
choose $\displaystyle \varepsilon \in \bigl(0, (1-(1-\delta)^p)^{\frac{q}{dq-dp}}\bigr).$
For $f(x):=|x|^{-d/q}$ and $k\in\mathbb{N}$, put $f_{\varepsilon,k}:=f\chi_{\{x:\varepsilon^{k+1}<|x|<\varepsilon^{k}\}}$.
Now write
\begin{align*}
f_1&:=(+1,+1,+1,+1):=f_{\varepsilon,3}+f_{\varepsilon,2}+f_{\varepsilon,1}+f_{\varepsilon,0},\\
f_2&:=(+1,+1,-1,-1):=f_{\varepsilon,3}+f_{\varepsilon,2}-f_{\varepsilon,1}-f_{\varepsilon,0},\\
f_3&:=(+1,-1,+1,-1):=f_{\varepsilon,3}-f_{\varepsilon,2}+f_{\varepsilon,1}-f_{\varepsilon,0}.
\end{align*}
Observe that $\|f_i\|_{\mathcal{M}^p_q}=\|f\chi_{\{x:\varepsilon^4<|x|<1\}}\|_{\mathcal{M}^p_q}$ for $i=1,2,3$, and that
\begin{align*}
3f_{\varepsilon,3} &\le |f_1+f_2+f_3| \le 3f\chi_{\{x:\varepsilon^4<|x|<1\}},\\
3f_{\varepsilon,2} &\le |f_1+f_2-f_3| \le 3f\chi_{\{x:\varepsilon^4<|x|<1\}},\\
3f_{\varepsilon,1} &\le |f_1-f_2+f_3| \le 3f\chi_{\{x:\varepsilon^4<|x|<1\}},\\
3f_{\varepsilon,0} &\le |f_1-f_2-f_3| \le 3f\chi_{\{x:\varepsilon^4<|x|<1\}}.
\end{align*}
By virtue of Corollary \ref{cor:140320-1}, we have
\[
3(1-\varepsilon^{d-dp/q})^{1/p}\|f\|_{\mathcal{M}^p_q} \le \|f_1\pm f_2 \pm f_3\|_{\mathcal{M}^p_q} \le
3\|f\chi_{\{x:\varepsilon^4<|x|<1\}}\|_{\mathcal{M}^p_q}
\]
for all choices of signs. For $i=1,2,3$, define $ \displaystyle F_i:=\frac{f_i}{\|f_i\|_{\mathcal{M}^p_q}}$.
Then, $\|F_i\|_{\mathcal{M}^p_q}=1$ for $i=1,2,3$, and
\[
\|F_1\pm F_2 \pm F_3\|_{\mathcal{M}^p_q}
=\frac{\|f_1\pm f_2 \pm f_3\|_{\mathcal{M}^p_q}}{\|f\chi_{\{x:\varepsilon^4<|x|<1\}}\|_{\mathcal{M}^p_q}}
\ge
\frac{3(1-\varepsilon^{d-dp/q})^{1/p}\|f\|_{\mathcal{M}^p_q}}{\|f\|_{\mathcal{M}^p_q}}
>3(1-\delta).
\]
This proves that $\mathcal{M}^p_q$ is not uniformly non-$\ell^1_3$.

In order to reveal the pattern, we shall now present the proof for $n=4$. With similar notations as above, we write
\begin{align*}
f_1&:=(+1,+1,+1,+1,+1,+1,+1,+1),\\
f_2&:=(+1,+1,+1,+1,-1,-1,-1,-1),\\
f_3&:=(+1,+1,-1,-1,+1,+1,-1,-1),\\
f_4&:=(+1,-1,+1,-1,+1,-1,+1,-1),
\end{align*}
where the $i$-th  term corresponds to the sign of $f_{\varepsilon,{8-i}}$ for $i=1,2, \ldots,8$.
%and the last term corresponds to the sign of $f_{\varepsilon,0}$.
Observe that $\|f_i\|_{\mathcal{M}^p_q}=\|f\chi_{\{x:\varepsilon^8<|x|<1\}}\|_{\mathcal{M}^p_q}$ for $i=1,\dots,4$, and that
\begin{align*}
4f_{\varepsilon,7} &\le |f_1+f_2+f_3+f_4| \le 4f\chi_{\{x:\varepsilon^8<|x|<1\}},\\
4f_{\varepsilon,6} &\le |f_1+f_2+f_3-f_4| \le 4f\chi_{\{x:\varepsilon^8<|x|<1\}},\\
4f_{\varepsilon,5} &\le |f_1+f_2-f_3+f_4| \le 4f\chi_{\{x:\varepsilon^8<|x|<1\}},\\
4f_{\varepsilon,4} &\le |f_1+f_2-f_3-f_4| \le 4f\chi_{\{x:\varepsilon^8<|x|<1\}},\\
4f_{\varepsilon,3} &\le |f_1-f_2+f_3+f_4| \le 4f\chi_{\{x:\varepsilon^8<|x|<1\}},\\
4f_{\varepsilon,2} &\le |f_1-f_2+f_3-f_4| \le 4f\chi_{\{x:\varepsilon^8<|x|<1\}},\\
4f_{\varepsilon,1} &\le |f_1-f_2-f_3+f_4| \le 4f\chi_{\{x:\varepsilon^8<|x|<1\}},\\
4f_{\varepsilon,0} &\le |f_1-f_2-f_3-f_4| \le 4f\chi_{\{x:\varepsilon^8<|x|<1\}}.
\end{align*}
Taking the Morrey norms, we get
\begin{align}%\label{eq:140320-1}
4(1-\varepsilon^{d-dp/q})^{1/p}\|f\|_{\mathcal{M}^p_q} \le \|f_1\pm f_2 \pm f_3 \pm f_4\|_{\mathcal{M}^p_q} \le
4\|f\chi_{\{x:\varepsilon^8<|x|<1\}}\|_{\mathcal{M}^p_q}
\end{align}
for all choices of signs. Taking $F_i:=\frac{f_i}{\|f_i\|_{\mathcal{M}^p_q}}$, we obtain
$\|F_i\|_{\mathcal{M}^p_q}=1$ for $i=1,\dots,4$. By our choice of $\varepsilon$ and the
fact that $\|f_i\|_{\mathcal{M}^p_q}=\|f\chi_{\{x:\varepsilon^8<|x|<1\}}\|_{\mathcal{M}^p_q}$
for $i=1,\dots,4$, we get
\[
\|F_1\pm F_2 \pm F_3 \pm F_4\|_{\mathcal{M}^p_q} > 4(1-\delta).
\]
Hence $\mathcal{M}^p_q$ is not uniformly non-$\ell^1_4$. Continuing the pattern, we see that
$\mathcal{M}^p_q$ is not uniformly non-$\ell^1_n$ for $n\ge2$.
\end{proof}

\section{$n$-th Von Neumann-Jordan Constant and $n$-th James Constant}

In this section, let $n\ge2$. For a Banach space $(X,\|\cdot\|_X)$, the $n$-th Von Neumann-Jordan constant
$C_{{\rm NJ}}^{(n)}(X)$ \cite{KTH} and the $n$-th James constant $C_{{\rm J}}^{(n)}(X)$ \cite{MNPZ} are defined by
\[
C_{{\rm NJ}}^{(n)}(X):=\sup \left\{\frac{\sum_{\pm} \|x_1\pm\cdots \pm x_n\|_X^2}{2^{n-1}\sum_{i=1}^n \|x_i\|_X^2}\,:\,x_i\not=0,
\ i=1,\dots,n\right\},
\]
and
\[
C_{{\rm J}}^{(n)}(X):=\sup \bigl\{\min \{\|x_1\pm \cdots \pm x_n\|_X:{\rm all~possible~choices~of~signs}\}:
x_i\in S_X,\ i=1,\dots,n\bigr\},
\]
respectively. In the definition of $C_{{\rm NJ}}^{(n)}(X)$, the sum $\sum_{\pm}$ is taken over all possible choices of signs.

We state some results about the two constants. The last one is specific for Morrey spaces.

\medskip

\begin{theorem}\label{NJ-constants}\cite{KTH}
For a Banach space $(X,\|\cdot\|_X)$ in general, we have $1\le C_{{\rm NJ}}^{(n)}(X) \le n$.
In particular, $C_{{\rm NJ}}^{(n)}(X)=1$ if and only if $X$ is a Hilbert space.
\end{theorem}

\medskip

\begin{theorem}\label{J-constant}\cite{MNPZ}
For a Banach space $(X,\|\cdot\|_X)$ in general, we have $1\le C_{{\rm J}}^{(n)}(X) \le n$.
If dim$(X)=\infty$, then $\sqrt{n}\le C_{{\rm J}}^{(n)}(X) \le n$.
For a Hilbert space $(X,\langle\cdot,\cdot\rangle_X)$, we have $C_{{\rm J}}^{(n)}(X)=\sqrt{n}$.
\end{theorem}

\medskip

\begin{theorem}\label{J-NJ-ineq}
For a Banach space $(X,\|\cdot\|_X)$ in general, we have
\[
[C_{{\rm J}}^{(n)}(X)]^2 \le nC_{{\rm NJ}}^{(n)}(X).
\]
\end{theorem}

\medskip

\begin{proof}
For every $x_i\in S_X,\ i=1,\dots,n$, let $m:=\min \{\|x_1\pm \cdots \pm x_n\|_X:{\rm all~possible~choices~of~signs}\}$.
Then, clearly $m\le \Bigl( \prod_\pm \|x_1\pm \cdots \pm x_n\|_X\Bigr)^{1/n}$, where the product is taken over all
possible choices of signs. Next, by the GM-QM inequality and the last inequality, we have
\begin{align*}
%\Bigl( \prod_\pm \|x_1\pm \cdots \pm x_n\|_X\Bigr)^{1/n}
m
&\le \left(\frac{\sum_{\pm} \|x_1\pm\cdots \pm x_n\|^2_X}{2^{n-1}}\right)^{1/2}\\
&=\left(n\cdot \frac{\sum_{\pm} \|x_1\pm\cdots \pm x_n\|^2_X}{2^{n-1}\sum_{i=1}^n \|x_i\|_X^2}\right)^{1/2}\\
&\le \bigl(nC_{{\rm NJ}(n)(X)}\bigr)^{1/2}.
\end{align*}
Taking the supremum over all $x_i\in S_X,\ i=1,\dots,n$, the desired inequality follows.
\end{proof}

\medskip

\begin{theorem}
For $1\le p<q<\infty$, we have $C_{{\rm J}}^{(n)}(\mathcal{M}^p_q)=C_{{\rm NJ}}^{(n)}(\mathcal{M}^p_q)=n$.
\end{theorem}

\medskip

\begin{proof}
It follows immediately from Theorem \ref{unl1n} that
$C_{{\rm J}}^{(n)}(\mathcal{M}^p_q)=n$. Combining this fact and Theorem \ref{J-NJ-ineq},
we get $C_{{\rm NJ}}^{(n)}(\mathcal{M}^p_q)\ge n$.
On the other hand, by Theorem \ref{NJ-constants}, we have
$C_{{\rm NJ}}^{(n)}(\mathcal{M}^p_q)\le n$. Thus, $C_{{\rm NJ}}^{(n)}(\mathcal{M}^p_q)=n$.
\end{proof}

\section{Concluding Remarks}

Before we end our paper, let us consider a Banach space $(X,\|\cdot\|_X)$ which is {\it uniformly $n$-convex},
that is, for every $\varepsilon\in(0,n)$ there exists $\delta\in(0,1)$ such that for every $x_1,\dots,x_n\in S_X$
with $\|x_1\pm\cdots\pm x_n\|_X$ $>\varepsilon$ for all choices of signs except for $\|x_1+\cdots+x_n\|_X$, we have
$\|x_1+\cdots+x_n\|_X\le n(1-\delta)$. This condition is stronger than the uniformly non-$\ell^1_n$ condition, as
we state in the following theorem.

\medskip

\begin{theorem}
If $X$ is uniformly $n$-convex, then $X$ is uniformly non-$\ell^1_n$.
\end{theorem}

\medskip

\begin{proof}
Assuming that $X$ is uniformly $n$-convex, we have to find a $\delta>0$ such that for every $x_1,\dots,x_n\in S_X$
with $\|x_1\pm\cdots\pm x_n\|_X>n(1-\delta)$ for all choices of signs except for $\|x_1+\cdots+x_n\|_X$, we have
$\|x_1+\cdots+x_n\|_X\le n(1-\delta)$. To do so, just take an $\varepsilon\in(0,n)$ and choose a corresponding
$\delta\in(0,1)$ such that for every $x_1,\dots,x_n\in S_X$ with $\|x_1\pm\cdots\pm x_n\|_X>\varepsilon$ for all
choices of signs except for $\|x_1+\cdots+x_n\|_X$, we have $\|x_1+\cdots+x_n\|_X\le n(1-\delta)$. Observe that if
$n(1-\delta)\ge\varepsilon$, then we are done. Otherwise, we choose $\delta_0\in(0,\delta)$ such that $n(1-\delta_0)
\ge\varepsilon$. This $\delta_0$ satisfies the uniformly non-$\ell^1_n$ condition.
\end{proof}

As a consequence of the above theorem and the fact that, for $1\le p<q<\infty$, the Morrey space $\mathcal{M}^p_q$ 
is not uniformly non-$\ell^1_n$, we conclude that $\mathcal{M}^p_q$ is not uniformly $n$-convex.

\bigskip

\noindent{\bf Acknowledgements}. This work is supported by P3MI-ITB 2020 Program. We thank
H.~Batkunde, Ifronika, and N.K.~Tumalun for useful discussions regarding the $n$-th James
constant and $n$-th Von Neumann-Jordan constant for Banach spaces.


\medskip

\begin{thebibliography}{999}

\bibitem{clarkson}
J.~A.~Clarkson,
The von Neumann--Jordan constant for the Lebesgue spaces,
{\em Ann. Math.} {\bf 38} (1937), No.~1, 114--115.

\bibitem{gao1990geometry}
J.~Gao and K.-S.~Lau,
On the geometry of spheres in normed linear spaces,
{\em J. Austral. Math. Soc.} {\bf 48} (1990), No.~1, 101--112.

\bibitem{GKSS}
H.~Gunawan, E.~Kikianty, Y.~Sawano, and C.~Schwanke,
Three geometric constants for Morrey spaces,
{\em Bull. Korean Math. Soc.} {\bf 56} (2019), No.~6, 1569--1575.

\bibitem{james}
R.~C.~James, Uniformly non-square Banach spaces,
{\em Ann. Math.} {\bf 80} (1964), No.~3, 542--550.

\bibitem{jimenez}
A.~Jim{\'e}nez-Melado, E.~Llorens-Fuster, and E.~Mazcun{\'a}n-Navarro,
The Dunkl--Williams constant, convexity, smoothness and normal structure,
{\em J. Math. Anal. Appl.} {\bf 342} (2008), No.~1, 298--310.

\bibitem{jordan1935inner}
P.~Jordan and J.~Von~Neumann,
On inner products in linear, metric spaces,
{\em Ann. Math.} {\bf 36} (1935), 719--723.

\bibitem{kato}
M.~Kato, L.~Maligranda, and Y.~Takahashi,
On James and Jordan--von Neumann constants and the normal structure
coefficient of Banach spaces,
{\em Studia Math.} {\bf 144} (2001), 275--295.

\bibitem{KTH}
M.~Kato, Y.~Takahashi, and K.~Hashimoto,
On $n$-th Von Neumann-Jordan constants for Banach spaces,
{\em Bull. Kyushu Inst. Tech.} No. 45 (1998), 25--33.

\bibitem{MNPZ}
L.~Maligandra, L.I.~Nikolova, L.-E.~Persson, T.~Zachariades,
On $n$-th James and Khintchine constants of Banach spaces,
{\em Math. Ineq. Appl.} {\bf 11}(1) (2007), 1--22.

\bibitem{MG}
A.~Muta'zili and H.~Gunawan,
On geometric constants for (small) Morrey spaces,
submitted.

\bibitem{sawano}
Y.~{Sawano},
A thought on generalized Morrey spaces,
{\em J. Indones. Math. Soc.} {\bf 25} (2019), 210--281.

\end{thebibliography}
\end{document}